\documentclass[conference]{IEEEtran}
\IEEEoverridecommandlockouts
\usepackage{tikz,pgfplots}
\usetikzlibrary{datavisualization}
\usetikzlibrary{datavisualization.formats.functions}
\usetikzlibrary{plotmarks}

\usepackage{cite}
\usepackage{amsmath,amssymb,amsfonts}
\usepackage{alphalph}
\usepackage{algorithmic}
\usepackage{graphicx}
\usepackage{textcomp}
\usepackage{xcolor}
\usepackage{eurosym}
\usepackage{multirow}
\usepackage{array}
\def\BibTeX{{\rm B\kern-.05em{\sc i\kern-.025em b}\kern-.08em
    T\kern-.1667em\lower.7ex\hbox{E}\kern-.125emX}}
\usepackage{bm}
\bmdefine{\blambda}{\lambda}
\bmdefine{\bmu}{\mu}

\usepackage{fancyhdr}
\fancypagestyle{specialfooter}{%
  \fancyhf{}
  
  \fancyfoot[C]{978-1-7281-1156-8/19/\$31.00 \copyright 2019 IEEE}
}

\begin{document}

\title{Robust Transmission Network Expansion Planning Problem Considering Storage Units\\
}

\author{\IEEEauthorblockN{\'Alvaro Garc\'ia-Cerezo}
\IEEEauthorblockA{\textit{Department of Electrical Engineering} \\
\textit{Universidad de Castilla-La Mancha}\\
Ciudad Real, Spain \\
Alvaro.Garcia29@alu.uclm.es}
\and
\IEEEauthorblockN{Luis Baringo}
\IEEEauthorblockA{\textit{Department of Electrical Engineering} \\
\textit{Universidad de Castilla-La Mancha}\\
Ciudad Real, Spain \\
Luis.Baringo@uclm.es}
\and
\IEEEauthorblockN{Raquel Garc\'ia-Bertrand}
\IEEEauthorblockA{\textit{Department of Electrical Engineering} \\
\textit{Universidad de Castilla-La Mancha}\\
Ciudad Real, Spain \\
Raquel.Garcia@uclm.es}

\thanks{

This work has been partially funded by the Ministry of Science, Innovation, and Universities of Spain under Projects RTI2018-096108-A-I00 and RTI2018-098703-B-I00 (FEDER, UE), the Ministry of Education and Professional Training of Spain under Grant 998142, and the Universidad de Castilla-La Mancha under Grant BI2018.
}

}

\maketitle

\begin{abstract}
This paper addresses the transmission network expansion planning problem considering storage units  under uncertain demand and generation capacity.
A two-stage adaptive robust optimization framework is adopted whereby short- and long-term  uncertainties are accounted for.
This work differs from previously reported solutions in an important aspect, namely, we include binary recourse variables to avoid the simultaneous charging and discharging of storage units once uncertainty is revealed.
Two-stage robust optimization with discrete recourse problems is a challenging task, so we propose using a nested column-and-constraint generation algorithm to solve the resulting problem. This algorithm guarantees convergence to the global optimum in a finite number of iterations.
The performance of the proposed algorithm is illustrated using the Garver's test system.
\end{abstract}


\section{Introduction}\label{section1}

\thispagestyle{specialfooter}

The Transmission Network Expansion Planning (TNEP) problem identifies the optimal reinforcements to be made in the transmission network of a power system.
This problem is generally solved by a central entity, such as the system operator, that determines the transmission network investment decisions that are optimal for the power system as a whole, e.g., those investment decisions that minimize investment and operating costs.
Since the pioneering work by Garver~\cite{Garver70}, the relevance of the TNEP problem has motivated the development of many models in the technical literature.

Over the last decade, the integration of renewable energy resources in power systems has given rise to the development of TNEP models that consider an uncertain environment in the decision-making process\cite{Conejo16}.
In this regard, two frameworks have been used to deal with this type of problem.
Firstly, stochastic programming \cite{Birge11}, which characterizes the uncertainty using a discrete set of realizations or scenarios.
This uncertainty framework is adequate if the probabilistic distribution of uncertain parameters is known; however, this task is not trivial.
On the other hand, the number of scenarios needed to obtain an accurate representation of the uncertainty is generally very large, which may result in tractability issues.
%
%
Secondly, robust optimization \cite{BenTal09}, which overcomes tractability issues associated with stochastic programming, but in contrast, solutions might be too conservative. 
Nevertheless, computational tractability is more restrictive.
Moreover, one of the aims of the TNEP problem is guaranteeing the supply of demands in all situations so that a conservative solution is not a drawback in this type of problems.
Thus the robust framework is preferred by scientist and researchers to deal with this type of problems. We also advocate and use a robust framework.
In particular, we propose a two-stage adaptive robust optimization approach. 

There are many references in the technical literature about the TNEP problem considering a robust optimization framework, such as, \cite{Jabr13,Roldan18}.
However, none of them consider storage units, and we consider that this is a significant aspect, especially since penetration of renewable units in power systems is expected to significantly increase in the next years. Thus, it is relevant to consider the impact of including storage units to store the renewable energy production during low-demand periods to be of use when is needed.

The TNEP problem considering storage units under a robust framework is analyzed in \cite{Dehghan16,Zhang18,Dvorkin18}.
Prevention of the simultaneous charge and discharge of storage units is modeled in \cite{Dehghan16} through first-stage binary variables. However, this assumption is unrealistic because the proper way to deal with these variables must be once uncertainty is revealed during the recourse problem, i.e., considering them as second-stage variables. In contrast, \cite{Zhang18} and \cite{Dvorkin18} do not include binary variables to deal with the issue of simultaneous charging and discharging of storage units. They basically assume that this is not relevant or significant.
As shown in \cite{Li16}, if the storage charging price at any bus with storage units is higher than or equal to the locational marginal price, simultaneous charging and discharging is likely to occur. For this reason, in systems with high penetration of renewable generating units where there are periods with very low or null marginal prices, it is necessary to prevent this undesirable situation.

The consideration of binary variables in the lower-level problem prevents the use of the traditional column-and-constraint generation algorithm to solve the two-stage adaptive robust optimization problem, because the two lowermost optimization levels cannot be recast as an equivalent single-level optimization problem. We propose solving the resulting model using the nested column-and-constraint generation algorithm described in \cite{Zhao12}, which to our knowledge is the only exact solution procedure to date for this kind of problems.

In summary, the main contribution of this paper is to propose an adaptive robust optimization approach with integer recourse variables for the TNEP problem, which avoids the simultaneous charging and discharging of storage units and makes the problem tractable even for medium-large size systems. In addition, we also consider the impact of both short- and long-term uncertainties in the decision-making process.
Short-term uncertainties in demand and renewable production are considered through a set of representative days, while long-term uncertainties in demand growth and capacity of generating units are represented using confidence bounds.

The rest of the paper is organized as follows. Section \ref{section2} provides the formulation of the problem. Section \ref{section3} explains the solution procedure. Results from a case study are reported in Section \ref{section4}. Finally, Section \ref{section5} concludes the paper with some relevant remarks.

\section{Compact Problem Formulation}\label{section2}

The TNEP problem is formulated using the following three-level adaptive robust optimization formulation in compact form:
\begin{subequations} \label{CPF1}
	\renewcommand{\theequation}{\theparentequation \alphalph{\value{equation}}}
	\allowdisplaybreaks
	\begin{align}
	&\hspace{-3cm}\min_{\textbf{x},c^{\rm{wc}}} \; \textbf{e}^{T}\textbf{x} + c^{\rm{wc}} \label{CPF1a}\hspace{1cm}\\
	&\hspace{-3cm}\text{subject to:}\nonumber\\
	&\hspace{-3cm}\textbf{K}\textbf{x} = \textbf{f} \label{CPF1b}\\
	&\hspace{-3cm}\textbf{L}\textbf{x} \leq \textbf{g} \label{CPF1c}\\
	&\hspace{-3cm}\textbf{x} \in \mathbb{Z}^{n} \label{CPF1d}
	\end{align}
\end{subequations}
\vspace{-0.7cm}
\begin{subequations} \label{CPF2}
	\renewcommand{\theequation}{\theparentequation \alphalph{\value{equation}}}
	\allowdisplaybreaks
	\begin{align}
	&\hspace{-1.3cm}c^{\rm{wc}} = \Bigg\lbrace \max_{\textbf{u},c}\;c& \label{CPF2a}\\
	&\hspace{-0.8cm}\text{subject to:}\nonumber\\
	&\hspace{-0.8cm}\textbf{u} \in \cal{U} \label{CPF2b}
	\end{align}
\end{subequations}
\vspace{-0.7cm}
\begin{subequations} \label{CPF3}
	\renewcommand{\theequation}{\theparentequation \alphalph{\value{equation}}}
	\allowdisplaybreaks
	\begin{align}
	&\hspace{-0.1cm}c = \Bigg[ \min_{\textbf{y},\textbf{z}} \; \textbf{b}^{T}\textbf{y}\label{CPF3a}\\
	&\hspace{0.4cm}\text{subject to:}\nonumber\\
	&\hspace{0.4cm}\textbf{A}\textbf{x}+\textbf{B}(\textbf{x})\textbf{y}+\textbf{D}\textbf{u} = \textbf{a} : \blambda \label{CPF3b}\\	
	&\hspace{0.4cm}\textbf{F}\textbf{x}+\textbf{G}\textbf{y}+\textbf{H}\textbf{u}+\textbf{I}(\textbf{x})\textbf{z} \geq \textbf{d} : \bmu \label{CPF3c}\\
	&\hspace{0.4cm}\textbf{z}\in \{0,1\}^{m}\Bigg] \Bigg\rbrace, \label{CPF3d}
	\end{align}
\end{subequations}
where $\textbf{x}$ is a vector representing first-stage expansion decision variables; $c$ and $c^{\rm{wc}}$ are variables representing operating and worst-case operating cost, respectively; $\textbf{u}$ is the vector of worst-case uncertainty realizations; $\textbf{y}$ and $\textbf{z}$ are second-stage variables (continuous and binary) representing operating decisions; and $\blambda$ and $\bmu$ are the lower-level dual variables.
In our particular problem, the aim of the vector of binary variables $\textbf{z}$ is to avoid the simultaneous charging and discharging of storage units.
Additionally, $\textbf{A}$, $\textbf{D}$, $\textbf{F}$, $\textbf{G}$, $\textbf{H}$, $\textbf{K}$, and $\textbf{L}$ are coefficient matrices; $\textbf{B}(\textbf{x})$ and $\textbf{I}(\textbf{x})$ are matrices whose elements depend on expansion variables; and $\textbf{a}$, $\textbf{b}$, $\textbf{d}$, $\textbf{e}$, $\textbf{f}$, and $\textbf{g}$ are coefficient vectors. Finally, $\cal{U}$ represents the uncertainty set that contains all possible materializations of the uncertain parameters, $n$ is the dimension of vector $\textbf{x}$, and $m$ is the dimension of vector $\textbf{z}$.

Problem \eqref{CPF1}-\eqref{CPF3} is a three-level model. The upper-level \eqref{CPF1} determines the expansion decisions (first-stage variables) minimizing investment and operating costs; the middle-level \eqref{CPF2} detects the worst-case realizations of uncertainty sources for the investment plan identified by the upper-level, i.e., in this level uncertainty is revealed; and the lower-level \eqref{CPF3} allows system operators to take a recourse action optimizing operating costs once first-stage variables and uncertainty are known.

Upper-level problem \eqref{CPF1} minimizes the total cost \eqref{CPF1a}, which includes both the investment and worst-case operating costs.
Constraints \eqref{CPF1b} define which transmission lines are initially built in the network.
Constraints \eqref{CPF1c} impose an investment budget and limit the number of storage facilities to be built.
Constraints \eqref{CPF1d} set out the integer nature of the vector of expansion variables $\textbf{x}$.

The middle-level problem \eqref{CPF2} determines the worst-case uncertainty realization, i.e., the uncertainty realization that leads to the largest operating cost \eqref{CPF2a} once the expansion decisions are made by the upper-level.
Constraints \eqref{CPF2b} impose that uncertain variables (in our case demands and generation capacities) are characterized by uncertainty set $\cal{U}$.
It should be noted that we use a cardinality-constrained uncertainty set as described in \cite{minguez2016robust}.
In this reference, the uncertainty sets are modeled through the uncertainty budget $\Gamma$, representing the maximum number of random parameters that may reach their lower or upper limits.
In our model, we consider uncertainty in demands, conventional, and wind-power generating units, therefore we deal with three uncertainty budgets $\Gamma^{\rm{D}}$, $\Gamma^{\rm{G}}$, and $\Gamma^{\rm{W}}$, respectively.

The lower-level problem \eqref{CPF3} identifies the minimum operating cost given upper- and middle-level decisions; i.e., problem \eqref{CPF3} identifies the system operating decisions that minimize the operating cost \eqref{CPF3a}.
Constraints \eqref{CPF3b} and \eqref{CPF3c} model the operating feasibility set, comprising the power balance at each bus, capacity of conventional and renewable generating units, power flow limits through each transmission line, limits on load shedding, operation of each energy storage, limits on the energy storage levels of each storage unit, no simultaneous charging and discharging of energy storage units, limits on the charging and discharging power of storage units, and definition of reference bus.
Finally, constraints \eqref{CPF3d} set out the binary nature of the vector of variables $\textbf{z}$.

\section{Solution Approach}\label{section3}

Problem \eqref{CPF1}-\eqref{CPF3} is an instance of mixed-integer three-level programming with lower-level binary variables.
The proposed solution approach is the nested column-and-constraint generation algorithm described in \cite{Zhao12}.
This method involves two loops.
The outer loop comprises the iterative solution of a master problem and a max-min subproblem with lower-level binary variables.
The solution of the max-min subproblem involves the iterative solution of two optimization problems, namely the inner-loop master problem and the inner-loop subproblem.
Next, the master problem and the subproblem are presented.

\subsection{Master Problem}

The master problem constitutes a relaxation for the original problem \eqref{CPF1}-\eqref{CPF3} where the second-stage problem is iteratively approximated by a set of valid operating constraints.
The master problem at iteration $j$ of the outer loop is formulated as the following mixed-integer problem:
\begin{subequations} \label{CPF4}
	\renewcommand{\theequation}{\theparentequation \alphalph{\value{equation}}}
	\allowdisplaybreaks
	\begin{align}
	&\min_{\textbf{x},\textbf{y}_{i},\textbf{z}_{i},\eta} \textbf{e}^{T}\textbf{x} + \eta \label{CPF4a}\\
	&\text{subject to:}\nonumber\\
	&\textbf{K}\textbf{x} = \textbf{f} \label{CPF4b}\\
	&\textbf{L}\textbf{x} \leq \textbf{g} \label{CPF4c}\\
	&\textbf{x} \in \mathbb{Z}^{n} \label{CPF4d}\\
	&\eta \geq \textbf{b}^{T}\textbf{y}_{i}; \; i = 1,\,...\,,\,j-1\label{CPF4e}\\
	&\textbf{A}\textbf{x}+\textbf{B}(\textbf{x})\textbf{y}_{i}+\textbf{D}\textbf{u}^{(i)} = \textbf{a}; \; i = 1,\,...\,,\,j-1\label{CPF4f}\\	
	&\textbf{F}\textbf{x}+\textbf{G}\textbf{y}_{i}+\textbf{H}\textbf{u}^{(i)}+\textbf{I}(\textbf{x})\textbf{z}_{i} \geq \textbf{d}; \; i = 1,\,...\,,\,j-1\label{CPF4g}\\
	&\textbf{z}_{i}\in \{0,1\}^{m}; \; i = 1,\,...\,,\,j-1,\label{CPF4h}
	\end{align}
\end{subequations}
where the additional vectors of decision variables $\textbf{y}_{i}$ and $\textbf{z}_{i}$, respectively corresponding to $\textbf{y}$ and $\textbf{z}$, are associated with the uncertain realizations identified by the subproblem at outer-loop iteration $i$ through $\textbf{u}^{(i)}$.

The objective function \eqref{CPF4a} is equivalent to \eqref{CPF1a} except for the term $\eta$, which is an approximation of the worst-case operating cost.
Constraints \eqref{CPF4b}-\eqref{CPF4d} are equivalent to \eqref{CPF1b}-\eqref{CPF1d}, respectively.
Constraints \eqref{CPF4e} represents a lower bound for $\eta$.
Constraints \eqref{CPF4f}-\eqref{CPF4h} correspond to lower-level constraints \eqref{CPF3b}-\eqref{CPF3d}, respectively.

Note that non-linear terms $\textbf{B}(\textbf{x})\textbf{y}_{i}$ in constraints \eqref{CPF4e} include the product of binary expansion-decision variables of matrix $\textbf{B}(\textbf{x})$ and continuous operating variables $\textbf{y}$, while non-linear terms $\textbf{I}(\textbf{x})\textbf{z}_{i}$ in constraints \eqref{CPF4f} include the product of integer expansion-decision variables of matrix $\textbf{I}(\textbf{x})$ and binary variables $\textbf{z}$.
Both terms can be linearized as shown in \cite{sioshansi2017optimization} so that problem \eqref{CPF4} is finally recast as a mixed-integer linear programming problem.

\subsection{Subproblem}

At each iteration $j$ of the outer loop, the subproblem determines the worst-case uncertainty realizations yielding the maximum value of the minimum operating cost for a given upper-level decision provided by the previous master problem \eqref{CPF4}.
The subproblem is a max-min problem comprising the two lowermost levels \eqref{CPF2}-\eqref{CPF3} parameterized in terms of given upper-level decision variables $\textbf{x}^{(j)}$.
Due to the presence of binary variables in the lower level, we propose solving such a max-min problem through an inner loop comprising the iterative solution of two optimization problems as described in \cite{Zhao12}, namely the inner-loop master problem and the inner-loop subproblem.

\subsubsection{Inner-loop master problem}
The inner-loop master problem at iteration $j$ of the outer loop and iteration $\ell$ of the inner loop is formulated as follow:
\begin{subequations} \label{CPF5}
	\renewcommand{\theequation}{\theparentequation \alphalph{\value{equation}}}
	\allowdisplaybreaks
	\begin{align}
	&c^{\rm{wc}} = \max_{\xi,\textbf{u},\blambda_{k},\bmu_{k}} \xi \label{CPF5a}\\
	&\text{subject to:}\nonumber\\
	&\textbf{u} \in \cal{U} \label{CPF5b}\\
	&\textbf{B}^{T}(\textbf{x}^{(j)})\blambda_{k}+\textbf{G}^{T}\bmu_{k} = \textbf{b} ; \; k = 1,\,...\,,\,\ell \label{CPF5c}\\
	&\bmu_{k}\geq \textbf{0}; \; k = 1,\,...\,,\,\ell\label{CPF5d}\\	
	&\xi \leq (\blambda_{k})^{T}\left(\textbf{a}-\textbf{A}\textbf{x}^{(j)}-\textbf{D}\textbf{u}\right) + (\bmu_{k})^{T}\Big(\textbf{d} \nonumber\\
	&\hspace{5mm}-\textbf{F}\textbf{x}^{(j)}-\textbf{H}\textbf{u}-\textbf{I}(\textbf{x}^{(j)})\textbf{z}^{(k)}\Big); \; k = 1,\,...\,,\,\ell,\label{CPF5e}
	\end{align}
\end{subequations}
where $\bmu_{k}$ and $\blambda_{k}$ are additional vectors of decision variables associated with the fixed values of $\textbf{z}^{(k)}$, which are obtained from the solution of the inner-loop subproblems previously solved in the current inner loop iteration indexed by $k$.

The objective function \eqref{CPF5a} is the approximation of the worst-case operating cost at iteration $j$ of the outer-loop, which is maximized subject to:
\begin{enumerate}
  \item Constraints \eqref{CPF5b} comprising the characterization of the uncertainty set.
  \item Constraints \eqref{CPF5c} and \eqref{CPF5d}, which correspond to lower-level dual feasibility constraints for fixed values of $\textbf{z}$.
  \item Constraints \eqref{CPF5e}, where the value of the lower-level dual objective function for fixed values of $\textbf{z}$ is considered as an upper bound for $\xi$.
\end{enumerate}

Note that non-linear terms $(\bmu_{k})^{T}\textbf{H}\textbf{u}$ and $(\blambda_{k})^{T}\textbf{D}\textbf{u}$ included in constraints \eqref{CPF5e} are the product of binary uncertainty variables $\textbf{u}$ and continuous dual variables $\blambda_{k}$ and $\bmu_{k}$, respectively.
Both terms can be linearized as shown in \cite{sioshansi2017optimization}.

\subsubsection{Inner-loop subproblem}
The inner-loop subproblem at iteration $j$ of the outer loop and iteration $\ell$ of the inner loop is formulated as follows:
\begin{subequations} \label{CPF6}
	\renewcommand{\theequation}{\theparentequation \alphalph{\value{equation}}}
	\allowdisplaybreaks
	\begin{align}
	&c^{(\ell)} = \min_{\textbf{y},\textbf{z}} \; \textbf{b}^{T}\textbf{y}\label{CPF6a}\\
	&\text{subject to:}\nonumber\\
	&\textbf{A}\textbf{x}^{(j)}+\textbf{B}(\textbf{x}^{(j)})\textbf{y}+\textbf{D}\textbf{u}^{(\ell)} = \textbf{a} \label{CPF6b}\\	
	&\textbf{F}\textbf{x}^{(j)}+\textbf{G}\textbf{y}+\textbf{H}\textbf{u}^{(\ell)}+\textbf{I}(\textbf{x}^{(j)})\textbf{z} \geq \textbf{d} \label{CPF6c}\\
	&\textbf{z}\in \{0,1\}^{m}. \label{CPF6d}
	\end{align}
\end{subequations}

The objective function \eqref{CPF6a} is identical to \eqref{CPF3a}, where $c^{(\ell)}$ corresponds to the operating cost of the inner-loop iteration $\ell$.
Constraints \eqref{CPF6b}-\eqref{CPF6d} are identical to constraints \eqref{CPF3b}-\eqref{CPF3d}, where the values of the expansion-decision variables $\textbf{x}^{(j)}$ are obtained by the master problem \eqref{CPF4} at outer-loop iteration $j$ and the values of uncertainty realizations vector $\textbf{u}^{(\ell)}$ are identified by the inner-loop master problem \eqref{CPF5} at inner-loop iteration $\ell$.

\subsection{Algorithm}

The proposed nested column-and-constraint generation algorithm works as follow:
\begin{itemize}
	\item[1)] Initialization of the outer loop.
	\begin{itemize}
		\item[\textbullet] Initialize the outer-loop lower and upper bounds: $LB^{o} = -\infty$ and $UB^{o} = +\infty$.
		\item[\textbullet] Set the outer-loop iteration counter $j$ to 1.
		\item[\textbullet] Set the vector of uncertain variables $\textbf{u}^{(j)}$ to their nominal values.
	\end{itemize}
	\item[2)] Master problem solution.
	\begin{itemize}
		\item[\textbullet] Solve problem \eqref{CPF4}. Obtain the optimal solution of the vector of expansion variables $\textbf{x}^{(j)}$.
		\item[\textbullet] Update $LB^{o}$ to the objective function of the master problem \eqref{CPF4}: $LB^{o} \leftarrow \textbf{e}^{T}(\textbf{x}^{(j)})+\eta^{(j)}$.
	\end{itemize}
	\item[3)] Subproblem solution.
	\begin{itemize}
		\item[3.1)]  Initialization of the inner loop.
		\begin{itemize}
			\item[\textbullet] Initialize the inner-loop lower and upper bounds: $LB^{i} = -\infty$ and $UB^{i} = +\infty$.
			\item[\textbullet] Set the inner-loop iteration counter  $\ell$ to 1.
			\item[\textbullet] Set the vector of uncertain variables $\textbf{u}^{(\ell)}$ to their nominal values.
		\end{itemize}
		\item[3.2)] Solution of the inner-loop subproblem.
		\begin{itemize}
			\item[\textbullet] Solve problem \eqref{CPF6} for given $\textbf{u}^{(\ell)}$ and $\textbf{x}^{(j)}$. Obtain the optimal solution of the vector of binary variables $\textbf{z}^{(\ell)}$.
			\item[\textbullet] Update $LB^{i}$ to the maximum between the current value of $LB^{i}$ and the objective function of problem \eqref{CPF6}: $LB^{i} \leftarrow \max\{LB^{i},c^{(\ell)}\}$.
		\end{itemize}
		\item[3.3)] Solution of the inner-loop master problem.
		\begin{itemize}
			\item[\textbullet] Solve problem \eqref{CPF5} for given $\textbf{z}^{(k)}$, $k = 1,\,...\,,\,\ell$ and $\textbf{x}^{(j)}$. Obtain the optimal solution of the uncertain variables vector $\textbf{u}^{\ell}$.
			\item[\textbullet] Set $\textbf{u}^{(\ell+1)}\leftarrow \textbf{u}$.
			\item[\textbullet] Update $UB^{i}$ to the objective function of problem \eqref{CPF5}: $UB^{i}\leftarrow c^{\rm{wc}}$.
		\end{itemize}
		\item[3.4)] Inner-loop convergence checking. If $UB^{i}-LB^{i}$ is lower than a predefined tolerance $\varepsilon^{i}$, stop and set $\textbf{u}^{(j+1)}\leftarrow \textbf{u}^{(\ell+1)}$. Otherwise, increase the inner-loop iteration counter: $\ell \leftarrow \ell+1$ and go to step 3.2.
	\end{itemize}
	\item[4)] Adjustment of the outer-loop upper bound. Update $UB^{o}$ to the minimum between the current value of $UB^{o}$ and the sum of $UB^{i}$ and annualized investment costs: $UB^{o} \leftarrow \min\{UB^{o},\textbf{e}^{T}(\textbf{x}^{(j)})+UB^{i}\}$.
	\item[5)] Outer-loop convergence checking. If $UB^{o}-LB^{o}$ is lower than a predefined tolerance $\varepsilon^{o}$, stop. Otherwise, increase the outer-loop iteration counter: $j \leftarrow j+1$ and go to step 2.
\end{itemize}

\section{Case Study}\label{section4}

This section provides the results of a case study.

\subsection{Data}

\begin{figure}
	\centering
	\includegraphics[scale=0.074]{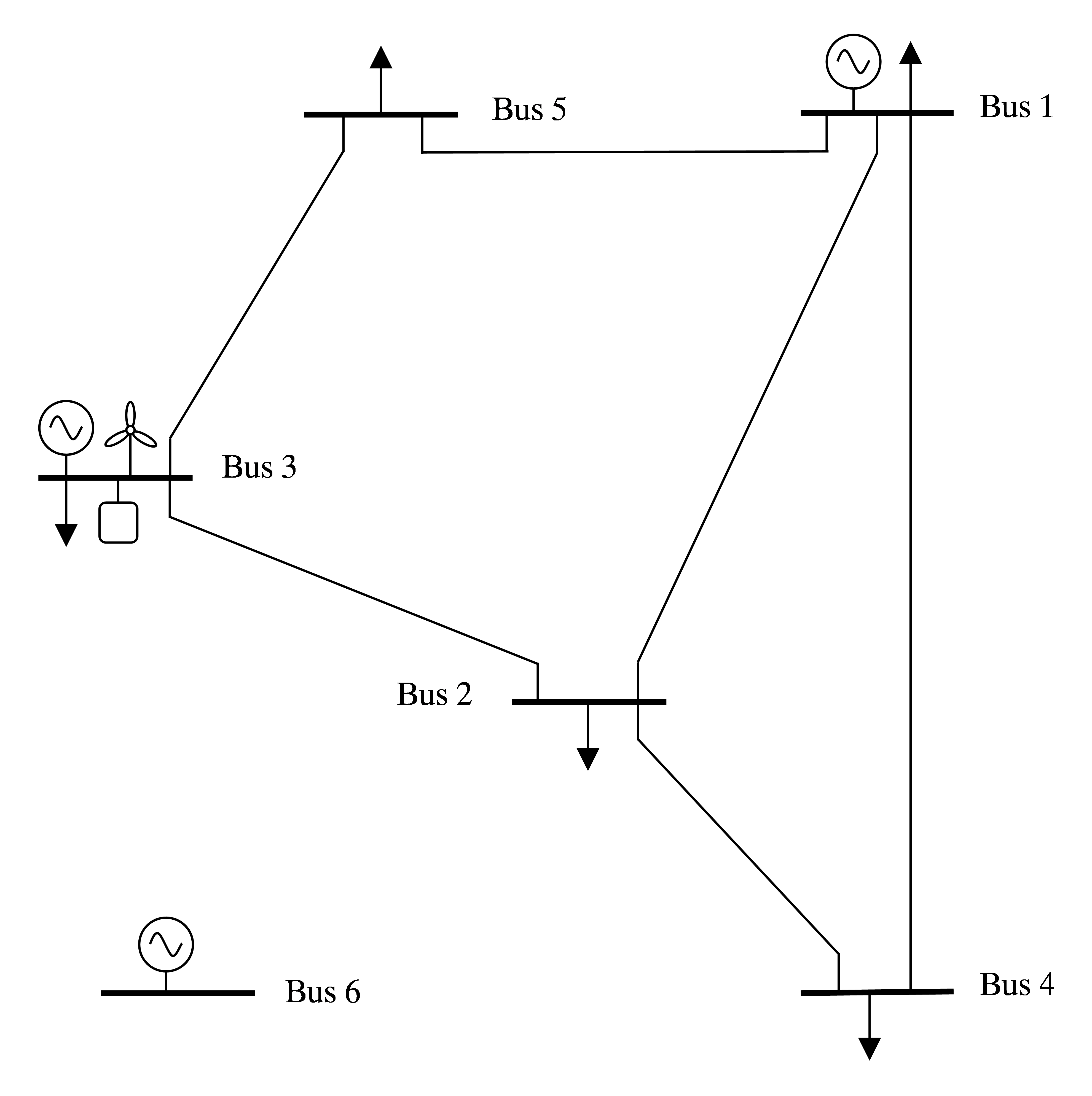}
	\caption{Garver's six-node test system.}
	\label{FigGarver}
\end{figure}

The proposed approach is analyzed using the Garver's six-node test system \cite{Garver70} depicted in Fig. \ref{FigGarver}.
This network comprises six buses, three conventional generating units, one wind-power generating unit, one existing storage unit, five demands, and six existing transmission lines.
Note that bus 6, which includes a conventional generating unit, is initially isolated. However, new transmission lines can be built connecting this bus with the rest of the network.
In particular, we consider that the number of existing and prospective transmission lines per corridor can be equal to three at most.
Additionally, we consider the possibility of building new storage units at bus 6.

Generating units, storage units, demand, and transmission line data are given in Tables~\ref{table_gen}, \ref{table_storage}, \ref{table_dem}, and \ref{table_lines}, respectively.

We consider that conventional and wind-power generating units can experience a maximum deviation of 50\% from their nominal capacities, while demand levels can increase up to a maximum of 20\%.
We consider an investment return period of 25 years and a discount rate of 10\%, which results in an annual amortization rate of 11\%.
We assume a total investment budget of \euro 60 million.
Furthermore, we consider that the convergence tolerances for the outer and the inner loops are $10^{-6}$.

Although  in this paper the formulation of the problem has focused on modeling the long-term uncertainty, special attention should be also given to short-term uncertainty, which plays a relevant role in the lower-level problem.
We consider that the short-term uncertainty has influence on the demand and wind-power production.
In order to model the short-term uncertainty, we apply a modified version of the K-means clustering technique described in \cite{baringo2013correlated} to the demand and wind historical data during 2016 in Texas \cite{energinet}.
We obtain 10 representative days of demand and wind-power production conditions, each one of them composed by 24 hourly data.

{
	\renewcommand{\arraystretch}{1.3}
	\begin{table}
		\centering
		\caption{Generating Unit Data}
		\label{table_gen}
		\begin{tabular}{l c c c} \hline
			& &Nominal&Operation\\
			Technology&Bus&capacity&cost\\
			& &(MW)&(\euro/MWh) \\ \hline
			Conventional&1&150&60\\
			Conventional&3&350&65\\	
			Conventional&6&500&70\\	
			Wind-power&3&400&0\\ \hline
		\end{tabular}
	\end{table}
}

{
	\renewcommand{\arraystretch}{1.6}
	\begin{table}
		\centering
		\caption{Storage Unit Data}
		\label{table_storage}
		\begin{tabular}{l c c} \cline{2-3}
			&Bus 3&Bus 6\\ \hline
			Maximum energy (MWh)&100&200\\
			Initial energy (MWh)&1&2\\
			Charging power capacity (MW)&20&40\\
			Discharging power capacity (MW)&20&40\\
			Charging efficienty (\%)&82&82\\
			Discharging efficienty (\%)&100&100\\
			Maximum number of&\multirow{2}{*}{-}&\multirow{2}{*}{3}\\
			units that can be built& & \\
			Investment cost ($10^3$\euro)&-&10,000\\ \hline
		\end{tabular}
	\end{table}
}

{
	\renewcommand{\arraystretch}{1.3}
	\begin{table}
		\centering
		\caption{Demand Data}
		\label{table_dem}
		\begin{tabular}{c c c} \hline
			\multirow{2}{*}{Bus}&Nominal&Load-shedding\\
			& level (MW)&cost (\euro/MWh)\\ \hline
			1&100&11,250\\
			2&300&11,500\\
			3&50&12,000\\
			4&200&11,000\\
			5&300&11,200\\ \hline
		\end{tabular}
	\end{table}
}

{
	\renewcommand{\arraystretch}{1.3}
	\begin{table}
		\centering
		\caption{Transmission Line Data}
		\label{table_lines}
		\begin{tabular}{c c c c c} \hline
			From&To&Reactance&Capacity&Investment\\
			Bus&Bus&(p.u.)&(MW)&cost ($10^3$\euro)\\ \hline
			1&2&0.40&100&7,723.20\\
			1&3&0.38&100&7,337.04\\
			1&4&0.60&80&11,584.80\\
			1&5&0.20&100&3,861.60\\
			1&6&0.68&70&13,129.44\\
			2&3&0.20&80&3,861.60\\
			2&4&0.40&100&7,723.20\\
			2&5&0.31&100&5,985.48\\
			2&6&0.30&100&5,792.40\\
			3&4&0.59&82&11,391.72\\
			3&5&0.20&70&3,861.60\\
			3&6&0.48&100&9,267.84\\
			4&5&0.63&75&12,164.04\\
			4&6&0.30&100&5,792.40\\
			5&6&0.61&78&11,777.88\\ \hline
		\end{tabular}
	\end{table}
}

Simulations were run using CPLEX 12.7.0.0 \cite{cplex} under GAMS 24.8.3 \cite{gamsguide} on a Gigabyte R280-A3C with 2 Intel Xeon E5-2698 at 2.3 GHz and 256 GB of RAM.

{
	\renewcommand{\arraystretch}{1.3}
	\begin{table}[h]
		\centering
		\caption{Case Study Results under Different Uncertainty Levels}
		\label{table_results}
		\begin{tabular}{c c c c c} \hline
			\multirow{3}{*}{Uncertainty level}&$\Gamma^{\rm{D}} = 0$&$\Gamma^{\rm{D}} = 2$&$\Gamma^{\rm{D}} = 3$&$\Gamma^{\rm{D}} = 5$\\
			$ $&$\Gamma^{\rm{G}} = 0$&$\Gamma^{\rm{G}} = 1$&$\Gamma^{\rm{G}} = 2$&$\Gamma^{\rm{G}} = 3$\\
			$ $&$\Gamma^{\rm{W}} \hspace{-0.06cm} = 0$&$\Gamma^{\rm{W}} \hspace{-0.06cm} = 0$&$\Gamma^{\rm{W}} \hspace{-0.06cm} = 1$&$\Gamma^{\rm{W}} \hspace{-0.06cm} = 1$\\ \hline
			 & &2-3 (x2)&1-5 (x1)& \\
			 &2-3 (x2)&2-5 (x1)&2-3 (x1)&2-3 (x1)\\
			New lines built&3-5 (x2)&3-5 (x2)&2-6 (x3)&2-6 (x3)\\
			 &4-6 (x2)&3-6 (x1)&3-5 (x2)&3-5 (x2)\\
			 & &4-6 (x3)&4-6 (x1)&\\
			 & &5-6 (x1)& & \\
			New storage units&\multirow{2}{*}{2}&\multirow{2}{*}{0}&\multirow{2}{*}{2}&\multirow{2}{*}{3}\\
			built at node 6& & & & \\
			Investment costs (M\euro)&47&60&59&59\\
			Total annual costs (M\euro)&94&162&1,497&2,624\\
			CPU (s)&28&49,255&16,549&636\\
			\# outer-loop iterations&1&5&4&2\\
			\# inner-loop iterations&1&2,$\,$2,$\,$2,$\,$2,$\,$3&2,$\,$3,$\,$2,$\,$2&11,$\,$2 \\ \hline
			
		\end{tabular}
	\end{table}
}

\subsection{Results}

We have solved the TNEP problem under different uncertainty levels, obtaining the expansion decisions, investment costs, total annual costs, computation times, and number of outer- and inner-loop iterations provided in Table \ref{table_results}. From this table we reach the conclusion that expansion decisions are clearly influenced by uncertainty levels.
In particular, building storage units is more attractive than building transmission lines for the case of maximum uncertainty level according to the results obtained.
Moreover, in all cases evaluated new lines are built connecting bus 6 with the rest of the network.

Regarding computational times, which correspond to the running time used until the outer-loop convergence is attained, we observe that the time needed is lower when the uncertainty budgets are null or equal to the number of demands, conventional, and wind-power generating units, respectively.
This is an expected result because in those cases the worst-case scenario is known in advance.
In addition, the number of outer-loop iterations needed to solve the problems follows the same behavior as the computational time.
Regarding the inner-loop problems, the number of iterations required to solve the subproblem is generally between one and three except for the case considering that all random parameters may reach their lower or upper limits, which requires 11 iterations during the first outer-loop iteration.

In order to demonstrate the importance of including binary variables in the model to avoid simultaneous charging and discharging of storage units, we have solved the problem without considering them. This problem has been solved using the traditional column-and-constraint generation algorithm, obtaining simultaneous charging and discharging in some cases, which highlights the importance of including these binary variables due to their potential impact on expansion decisions.
%


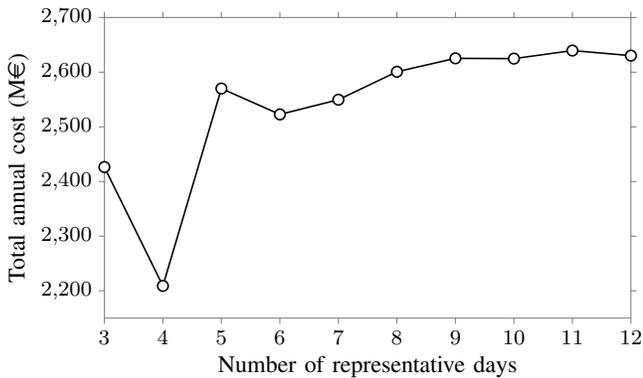
\begin{figure}[htb]
\begin{center}
\begin{tikzpicture}
\datavisualization [scientific axes,
x axis={attribute=t, include value={3}, label={Number of representative days}, ticks={step=1}, length=7cm,},
y axis={attribute=n, include value={2150,2700}, label={Total annual cost (M\euro)}, ticks={step=100}, length=4cm},
visualize as line=my data,
my data={style={mark=*,mark options={fill=white}}}]
data {
t, n
3, 2426.418
4, 2208.988
5, 2570.067
6, 2522.708
7, 2549.589
8, 2600.759
9, 2625.312
10, 2624.784
11, 2639.560
12, 2630.243
};
\end{tikzpicture}
	\caption{Total annual cost considering different number of representative days.}\label{FigNumberOfRP}
\end{center}
\end{figure}

Finally, note that the solution of the problem would be different depending on the number of representative days considered.
This fact has been analyzed using different numbers of representative days in order to find when the total annual cost of the TNEP problem becomes stable.
We have considered the uncertainty budgets $\Gamma^{\rm{D}}=5$, $\Gamma^{\rm{G}}=3$, and $\Gamma^{\rm{W}}=1$, obtaining the results depicted in Fig.~\ref{FigNumberOfRP}.
It can be seen that the total annual cost presents a stable value selecting a number of representative days greater than 8. For that reason, we have chosen 10 representative days to model the short-term uncertainty.

\section{Conclusions}\label{section5}

This paper proposes a new model to solve the TNEP problem considering storage units and avoiding the simultaneous charging and discharging of those units.
To that end, we adopt a two-stage adaptive robust optimization framework and include binary variables in the second stage variable set. The resulting three-level formulation forces us to use a nested column-and-constraint generation algorithm, whose description is also given.
Finally, we have analyzed the influence of the uncertainty level in expansion decisions, investment and operating costs, computational times, and number of outer- and inner-loop iterations on a case study using the Garver's six-node test system.

Given the theoretical framework and the results of the case study, the conclusions below are in order:
\begin{itemize}
	\item[1)] It is necessary to use binary variables in the third-level to avoid the simultaneous charging and discharging of storage units, which prevents the use of traditional column-and-constraint generation algorithms.
	\item[2)] Disregarding the prevention of simultaneous charging and discharging of storage units may have a great impact on expansion decisions, especially in systems with high renewable penetration.
    \item[3)] There is no clear pattern among uncertainty levels, expansion decisions, and the computational burden of the problem.
\end{itemize}

Future work comprises a more detailed description of the effect of simultaneous charging and discharging of storage units in larger systems, as well as the inclusion of transmission line contingencies in the formulation.

\end{document}